\documentclass[12pt]{amsart}
\usepackage{amscd}
\def\NZQ{\Bbb}               
\def\NN{{\NZQ N}}

\def\frk{\frak}               

\def\Phi{{\frk n}}
\def\Phi{{\frk N}}
\def\opn#1#2{\def#1{\operatorname{#2}}} 
\opn\chara{char} \opn\length{\ell} \opn\pd{pd} \opn\rk{rk}
\opn\projdim{proj\,dim} \opn\injdim{inj\,dim} \opn\rank{rank}
\opn\depth{depth} \opn\grade{grade} \opn\height{height}
\opn\embdim{emb\,dim} \opn\codim{codim}

\opn\Tr{Tr} \opn\bigrank{big\,rank}
\opn\superheight{superheight}\opn\lcm{lcm}
\opn\trdeg{tr\,deg}
\opn\reg{reg} \opn\lreg{lreg} \opn\ini{in} \opn\lpd{lpd}
\opn\size{size}
\opn\div{div} \opn\Div{Div} \opn\cl{cl} \opn\Cl{Cl}
\opn\Spec{Spec} \opn\Supp{Supp} \opn\supp{supp} \opn\Sing{Sing}
\opn\Ass{Ass} \opn\Min{Min}
\opn\Ann{Ann} \opn\Rad{Rad} \opn\Soc{Soc}
\opn\Im{Im} \opn\Ker{Ker} \opn\Coker{Coker} \opn\Am{Am}
\opn\Hom{Hom} \opn\Tor{Tor} \opn\Ext{Ext} \opn\End{End}
\opn\Aut{Aut} \opn\id{id}

\opn\nat{nat}
\opn\pff{pf}
\opn\Pf{Pf} \opn\GL{GL} \opn\SL{SL} \opn\mod{mod} \opn\ord{ord}
\opn\Gin{Gin} \opn\Hilb{Hilb}
\opn\aff{aff} \opn\con{conv} \opn\relint{relint} \opn\st{st}
\opn\lk{lk} \opn\cn{cn} \opn\core{core} \opn\vol{vol}
\opn\link{link} \opn\star{star}
\opn\gr{gr}

\def\pot#1#2{#1[\kern-0.28ex[#2]\kern-0.28ex]}

\opn\dirlim{\underrightarrow{\lim}}
\opn\inivlim{\underleftarrow{\lim}}
\def\Implies{\ifmmode\Longrightarrow \else
        \unskip${}\Longrightarrow{}$\ignorespaces\fi}
\def\implies{\ifmmode\Rightarrow \else
        \unskip${}\Rightarrow{}$\ignorespaces\fi}
\def\iff{\ifmmode\Longleftrightarrow \else
        \unskip${}\Longleftrightarrow{}$\ignorespaces\fi}

\let\:=\colon
\newtheorem{Theorem}{Theorem}[section]
\newtheorem{Lemma}[Theorem]{Lemma}
\newtheorem{Corollary}[Theorem]{Corollary}

\let\epsilon\varepsilon
\let\phi=\varphi
\let\kappa=\varkappa
\def\qed{\ifhmode\textqed\fi
      \ifmmode\ifinner\quad\qedsymbol\else\dispqed\fi\fi}
\def\textqed{\unskip\nobreak\penalty50
       \hskip2em\hbox{}\nobreak\hfil\qedsymbol
       \parfillskip=0pt \finalhyphendemerits=0}
\def\dispqed{\rlap{\qquad\qedsymbol}}

\opn\dis{dis}
\def\pnt{{\raise0.5mm\hbox{\large\bf.}}}

\opn\Lex{Lex}

\begin{document}

\title{Stanley Conjecture in small embedding dimension}

\author{Imran Anwar, Dorin Popescu}

\thanks{The authors are highly grateful to the School of Mathematical Sciences,
  GC University, Lahore, Pakistan in supporting and facilitating this
  research. The second author was supported by  CNCSIS and the Contract 2-CEX06-11-20/2006
 of the Romanian Ministery of Education and Research and the Higher Education Comission of Pakistan. }

\address{Imran Anwar, School of Mathematical Sciences, 68-B New Muslim Town,
    Lahore,Pakistan.}\email{iimrananwar@gmail.com}
\address{Dorin Popescu, Institute of Mathematics "Simion Stoilow",
University of Bucharest, P.O.Box 1-764, Bucharest 014700, Romania}
\email{dorin.popescu@imar.ro} \maketitle

\begin{abstract}
We show that Stanley's conjecture holds for a polynomial ring over a
field in four variables. In the case of polynomial ring in five
variables, we prove that the monomial ideals with all associated
primes of height two, are Stanley ideals.
 \vskip 0.4 true cm
 \noindent
\end{abstract}

\section{Introduction}

Let $S=K[x_1,x_2,...,x_n]$ be a polynomial ring in $n$ variables
over a field $K$ and $I\subset S$  a monomial ideal. In this paper a
{\em prime filtration} of $I$ is  assumed to be a monomial prime
filtration, that is a monomial filtration
\begin{center}
$\mathcal{F}$:\,\,\,\,\,\,\,\,\,\,\,\,\,\,\,\,\,\,\,\,\,\,\,\,$I=I_0\subset
                  I_1\subset .\,.\,.\subset I_r=S$
\end{center}
with   $I_j/I_{j-1}\cong S/P_j(-a_j)$ for some monomial prime ideals
$P_j$ of $S$, $a_j\in \NN^n$ and $j=1,2,...,r$. Set $\Supp(
\mathcal{F})=\{P_1,\ldots,P_r\}$. After \cite{HP}, the prime
filtration $\mathcal{F}$ is called {\em pretty clean}, if for all
$i<j$ for which $P_i\subseteq P_j$, it follows that $P_i=P_j$. The
monomial ideal $I$ is called {\em pretty clean}, if it has a pretty
clean filtration.

 Let $I\subset S$ be a monomial ideal, any
decomposition of $S/I$ as a direct sum of $K$-vector spaces of the
form $uK[Z]$ where $u$ is a monomial in $S$, and $Z\subseteq \{
x_1,x_2,...,x_n\} $ is called {\em Stanley decomposition}. Stanley
\cite{St} conjectured that there always exists a Stanley
decomposition
\begin{center}
                  $$S/I = \bigoplus_{i=1}^r u_iK[Z_i]$$
\end{center}
such that $|Z_i|\geq \depth (S/I)$ for all $i,1\leq i\leq r$. If
this is the case, we call $I$ a {\em Stanley ideal}. Sometimes
Stanley decompositions of $S/I$ arise from prime filtrations. In
fact, if $\mathcal{F}$ is a prime filtration of $S/I$ with factors
$(S/P_i)(-a_i)$ for $i=1,2,...,r$ then set $u_i=\prod_{j=1}^n x_j$
and $Z_i = \{x_j^{a_{ij}}: x_j\not \in P_i  \} $ and we have
$$S/I = \bigoplus_{i=1}^r u_iK[Z_i]$$
If $\mathcal{F}$ is a pretty clean filtration of $S/I$, then by
\cite[Corollary 3.4]{HP}
$$\Ass(S/I)=\Supp(\mathcal{F})$$ The converse is not always true
(see \cite[Example 4.4]{Ja}). However not all Stanley decompositions
arise from the prime filtrations (see \cite[Example 3.8]{MS}). A
prime filtration $\mathcal{F}$ is a {\em Stanley filtration} if the
Stanley decomposition arising from $\mathcal{F}$ satisfies the
Stanley conjecture. In \cite[Proposition 2.2]{Ja} it is  shown  that
all prime filtrations $\mathcal{F}$ for which
$\Ass(S/I)=\Supp(\mathcal{F})$ are Stanley filtrations, in
particular   all  monomial ideals $I\subset S$ for which $S/I$ is
pretty clean, are Stanley (see \cite[Theorem 6.5]{HP}). In case
$n=3$, for any monomial ideal $I\subset S$ we have $S/I$  pretty
clean by \cite[Theorem 1.10]{Ja} and so $I$ is Stanley. This result
was first obtained by different methods in \cite{Ap}. Recently,
Herzog, Soleyman Jahan, Yassemi \cite[Proposition 1.4]{HJY} showed
that if $I$ is a monomial ideal of $S$ (for any $n$) such that $S/I$
is Cohen-Macaulay of codimension two then $I$ is a Stanley ideal.

It is the purpose of our note to describe the Stanley ideals of the
polynomial ring $S=K[x_1,x_2,...,x_n]$, $n\leq 5 $. If $n=4$ we show
that all monomial ideals of $S$ are Stanley (see Theorem
\ref{main}). This extends \cite[Corollary 1.4]{AP}, which says that
a sequentially Cohen-Macaulay monomial ideal $I\subset
K[x_1,\ldots,x_4]$ is Stanley. If $n=5$ we show that all monomial
ideals $I\subset S$ having all the associated prime ideals of height
2 are Stanley ideals (see Corollary \ref{cor}).

\section{Stanley's Conjecture in small embedding dimension}

We start with a very elementary Lemma.
\begin{Lemma}\label{1}
Let $S=K[x_1,x_2,...,x_n]$, $T=K[x_1,x_2,...,x_r]$ for some $1\leq
r\leq n$ and $\mathcal{J}\subset T$ a monomial ideal. Then
$T/\mathcal{J}$ is pretty clean if and only if $S/\mathcal{J}S$ is
pretty clean.
\end{Lemma}

The following lemma is a key result in this note.

\begin{Lemma}\label{2}
Let $S=K[x_1,x_2,...,x_n]$, $n\leq 5 $ be a polynomial ring and
$I\subset S$  a monomial ideal having all the associated primes of
height $2$. Then there exists a prime filtration $\mathcal{F}$ of
$S/I$ such that $ht(P)\leq 3$ for all $P\in \Supp(\mathcal{F})$.
\end{Lemma}
\begin{proof}
We  use induction on $\mathbf{s}(I)$, where $\mathbf{s}(I)$ denotes
the number of irreducible monomial ideals appearing in the unique
decomposition of $I$ as an intersection of irreducible monomial
ideals (see \cite[Theorem 5.1.17]{Vi}), let us say $$ I =
\bigcap_{i=1}^{\mathbf{s}(I)} Q_i ,$$
 where $Q_i$'s are irreducible monomial ideals of codimension $2$.
If $\mathbf{s}(I)=1$, then the result follows because $S/I$ is
clean. If $\mathbf{s}(I)\geq 1$, then set $$ \mathcal{J}=
\bigcap_{i=2}^{\mathbf{s}(I)} Q_i\,\,\,.$$ Therefore $I =
\mathcal{J}\cap Q_1 $. We may suppose $Q_1=(x_1 ^{d_1},x_2^{d_2})$
after renumbering of variables, with $d_1$ the largest power of
$x_1$ in $\bigcup_{i=1}^{\mathbf{s}(I)} G(Q_i)$, where $G(Q_i)$  is
the set of minimal monomial generators of $I$. We claim that
$$\mathcal{F}_0=I\,\subset\,\,\mathcal{F}_1=(I,x_2^{d_2})\,\,\,{\subset\,\,\mathcal{F}_2=(x_1^{d_1},x_2^{d_2})}\,\subset\mathcal{F}_3= S $$
is a filtration of $S/I$, which will give the desired filtration by
refinating. Note that $\mathcal{F}_3 /\mathcal{F}_2 =
S/(x_1^{d_1},x_2^{d_2})$ is a clean module, so $\mathcal{F}_3
/\mathcal{F}_2 $ has a prime filtration involving only the prime
$(x_1,x_2)$.

 Now for $ \mathcal{F}_2
/\mathcal{F}_1=(x_1^{d_1},x_2^{d_2})/(I,x_2^{d_2})\cong
S/((I,x_2^{d_2}):x_1^{d_1})$ we have
$$E:=((I,x_2^{d_2}):x_1^{d_1})=
(I:x_1^{d_1},x_2^{d_2})=(\mathcal{J}:x_1^{d_1},x_2^{d_2})$$ and we
get
$$E=\bigcap_{i=2}^{\mathbf{s}(I)} ((Q_i : x_1^{d_1}),x_2^{d_2}).$$
Set $T=K[x_2,\ldots,x_{n}]$. Since $U_i:=((Q_i :
x_1^{d_1}),x_2^{d_2})$ is either $S$, or an irreducible ideal of
height $2$ or $3$ in the variables $x_2,\ldots,x_n$, we note that
$E=WS$ for a monomial ideal $W\subset T$ with all associated prime
ideals of dimension $n-2$ or $n-3$. If $n=4$ then $\dim T=3$ and
$T/W$ is pretty clean by \cite[Theorem 1.10]{Ja} and so $S/E$ is
pretty clean by Lemma \ref{1}. If $n=5$ then set $G:=\bigcap_{i=2,\
ht(U_i)=2}^{\mathbf{s}(I)}\ U_i$ and consider the filtration
$W\subset G\subset T$ (this is the dimension filtration of
\cite{Sc}). As $\mathbf{s}(G)< \mathbf{s}(I)$ we get by induction
hypothesis a prime filtration of $S/G$ involving just prime of
height $\leq 3$. Since $\Ass(G/W)$ contains just prime ideals of
height 3 we get $G/W$ clean by
 \cite[Corollary 2.2]{Po}. So we get a prime filtration of $T/W$
and by extension of $S/E$, involving only prime ideals of height
$\leq 3$. Therefore in both cases, $\mathcal{F}_2
/\mathcal{F}_1\cong S/E$ has a prime filtration involving only prime
ideals of height at most $3$.

 Finally,  $\mathcal{F}_1
/\mathcal{F}_0=(I,{x_2}^{d_2})/I\cong S/(I:x_{2}^{d_2})$ , and we
have
$$(I:x_2^{d_2})= (\mathcal{J}:x_2^{d_2})=\bigcap_{i=2}^{\mathbf{s}(I)} (Q_i : x_2^{d_2}),$$
where $(Q_i : x_2^{d_2})$ is either $S$, or irreducible  of  height
2. Since $ \mathbf{s}(I:x_2^{d_2})\leq
\mathbf{s}(I)-1<\mathbf{s}(I)$, we get by induction hypothesis that
$S/(I:x_2^{d_2})$ (and so $\mathcal{F}_1 /\mathcal{F}_0$ ) has a
prime filtration with prime ideals of height at most $3$. Then by
gluing together all these prime filtrations we get the desired one.
\end{proof}

\begin{Corollary}\label{cor}
Let $S=K[x_1,x_2,...,x_n]$, $n\leq 5 $ be a polynomial ring and
$I\subset S$  a monomial ideal having all the associated prime
ideals of height $2$. Then I is a Stanley ideal.
\end{Corollary}

\begin{proof}
If $S/I$ is Cohen-Macaulay then $I$ is a Stanley ideal by
\cite[Proposition 1.4]{HJY}. Now if $S/I$ is not Cohen-Macaulay
 then $\depth(S/I)\leq n-3$ because  $\dim(S/I)= n-2$. Let $\mathcal{F}$ be
the filtration given by  Lemma \ref{2}. Then all the associated
primes $P$ of $\Supp(\mathcal{F})$ satisfy the
condition,$$\dim(S/P)\geq n-3\geq \depth(S/I) \,\,.$$ Thus $I$ is a
Stanley ideal.
\end{proof}

\begin{Theorem}\label{main}
Any monomial ideal $I$ of $S=K[x_1,x_2,x_3,x_4]$ is a Stanley ideal.
\end{Theorem}
\begin{proof}
Let $\mathbf{s}(I)$ be the number of irreducible monomial ideals
appearing in the unique decomposition of $I$ as an intersection of
irreducible monomial
ideals,\\
let us say $$ I = \bigcap_{i=1}^4\bigcap_{j=1}^{\mathbf{s}(Q_i)}
Q_{ij} \,\,\,\  {\mbox{and}} \,\,\,
Q_i=\bigcap_{j=1}^{\mathbf{s}(Q_i)} Q_{ij} \  , \ \
\mathbf{s}(I)=\Sigma_{i=1}^4 {\mathbf{s}(Q_i)} ,$$ where $Q_{ij}$'s
are irreducible monomial ideals of height $i$. If $ht(I)=t$, $1\leq
t\leq 4$ then $Q_k=S$ and $\mathbf{s}(Q_k)=0$ for all $1< k< t$.
After Schenzel \cite{Sc}, the dimension filtration of I will be
$$\mathcal{F}_0=I\,\subset\,\,\mathcal{F}_1=
Q_1 \cap Q_2 \cap Q_3\,\,\,\subset\,\,\mathcal{F}_2= Q_1 \cap Q_2\,\subset\mathcal{F}_3=
 Q_1\subset\mathcal{F}_4= S $$
Now consider $\mathcal{F}_4 /\mathcal{F}_3\cong S/Q_1$, where
$Q_1=\bigcap_{j=1}^{\mathbf{s}(Q_1)}Q_{1j}$ with $Q_{1j}$'s
principal ideals. Therefore $Q_1=(u)$ for a monomial $u$ in $S$ (it
is factorial ring) and so  $S/Q_1$ is pretty clean (see e.g. the
proof of \cite[Lemma 1.9]{Ja}). Hence $\mathcal{F}_3$ is a Stanley
ideal.

 Now take
$\mathcal{F}_3 /\mathcal{F}_2\cong Q_1/(Q_1\cap Q_2)\cong S/(Q_2 :
u)$, where
$$(Q_2 : u)=\bigcap_{j=1}^{\mathbf{s}(Q_2)}(Q_{2j}:u)$$ Also
$(Q_{2j}:u)$ is either $S$ or of  height 2  for all $j$. By
Corollary \ref{cor},  $(Q_{2}:u)$ is a Stanley ideal and so is
$\mathcal{F}_2 $ because the clean filtrations of $S/(u)$ involve
only prime ideals of $\depth \ 3$. If $\mathbf{s}(Q_3)=
\mathbf{s}(Q_4)=0$ then we are done. If $\mathbf{s}(Q_3)\neq 0$ then
$\depth(S/I)\leq 1$ and $\Ass(\mathcal{F}_2/\mathcal{F}_1)$ contains
only prime ideals of height $3$. Hence $\mathcal{F}_2/\mathcal{F}_1$
is pretty clean by \cite [Corollary 2.2]{Po} and so
$\mathcal{F}_2/\mathcal{F}_1$ has a prime filtration involving only
prime ideals of height $3$. But as above $S/\mathcal{F}_2$ has a
prime filtration with prime ideals of height $\leq 3$. Gluing
together these two prime filtrations we get a prime filtration with
prime ideals $P$ such that $\dim (S/P)\geq 1\geq \depth(S/I)$. So
$\mathcal{F}_1$ is a Stanley ideal. If $\mathbf{s}(Q_4)\neq 0$ then
$\depth(S/I)=0$ and every prime filtration gives a Stanley
filtration.
\end{proof}

\end{document}